%% file: dxi.tex
\documentclass[11pt,eqno,sumlimits]{amsart}

\usepackage{amssymb, amscd, amsmath, epsfig, mathtools}
\usepackage[utf8]{inputenc}
\usepackage{amsthm}
\usepackage{enumerate}
\usepackage{xcolor}
\usepackage{scalerel}
\usepackage{soul}
\usepackage{tikz-cd}
\usepackage{esint}
\usepackage{slashed}
\usepackage[textsize=tiny]{todonotes}
\usepackage{pinlabel} % for labelling the figure
\usepackage{hyperref}
\usepackage{url}

\usepackage[margin=1.0in]{geometry}

\newcommand{\Z}{\mathbb Z}

\input{math_commands.tex}

\usepackage{hyperref}
\usepackage{url}
\usepackage{booktabs} 
\usepackage{xcolor}

\usepackage{times}
\usepackage{microtype}
\usepackage{graphicx}
\usepackage{subfigure}

\usepackage{tikz}    
\usepackage{tikz-3dplot}  
\usepackage{tikz-cd}
\usepackage[a]{esvect}
\usepackage{comment}
\usepackage{algorithm}
\usepackage{algpseudocode}
\usepackage{caption}

\title{A note on reproducing kernels for Sobolev spaces}

\author{
Steven Rosenberg}
\address{Department of Mathematics and Statistics,
Boston University,
Boston, MA 02215}
\email{sr@math.bu.edu}

\newtheorem{theorem}{Theorem}[section]
\newtheorem{lemma}{Lemma}[section]

\newtheorem{remark}[theorem]{Remark}

\def\calF{\mathcal{F}}

\def\R{\mathbb{R}}

\def\Z{\mathbb{Z}}

\def\calH{{\mathcal H}}

\newcommand{\dxix}{d_x(y)}

\def\Z{\mathbb Z}

\newcommand{\dxi}{d_{x}}

\newcommand{\vol}{{\rm vol}}

\begin{document}

\maketitle

\begin{abstract}
In this note, we compute the reproducing kernel for the RKHS of functions on $\R^n$ in a sufficiently high Sobolev norm.
\end{abstract}

\section{Introduction}An RKHS is a Hilbert space of complex-valued functions $\calH \subset \{f:M\to \R\}$ for some topological space $M$ with the property that the evaluation/delta functions $\delta_x(f) = f(x)$ are continuous for all $x\in M.$  Thus there exists $d_x\in \calH$ such that
$\delta_x(f) = \langle d_x, f\rangle_\calH$. For $M$ an $n$-dimensional closed manifold, the Sobolev space $H_s(M)$ is an RKHS for integers $s> {\rm dim}(M)/2.$  In this note, we explicitly compute $d_x$ for $M= \R^n$.
 This produces $d_x$ on a general manifold via partition of unity, but the explicitness is lost.

\section{The Computation}

Let $N(x,\sigma)(x) $ be the multivariate normal distribution with mean $x$ and variance $\sigma = \sigma\cdot {\rm Id}.$  For $\sigma\approx 0$, 
$N(x,\sigma)$ acts like a delta function: for $f\in H_s(\R^n),$ %, s > (n/2) +1,$
$$\lim_{\sigma\to 0} \int_{\R^n} N(x,\sigma)(x) f(x) dx = f(x).
$$
This uses that $f$ is continuous, so the delta function is a continuous functional on $H_s$ by the Sobolev embedding theorem. Hence, $f$ is in 
$H_{-s}$, 
and $\lim_{\sigma\to 0} N(x,\sigma) = \delta_{x}$ in $H_{-s}(\R^n).$ 
Therefore
\begin{align*}\lim_{\sigma\to 0} \int_{\R^n} N(x,\sigma) f(x) dx &= f(x) = \langle \dxi , f\rangle_s = \int_{\R^n}
\widehat{\dxi }(\xi)\hat f(\xi) (1+ |\xi|^2)^{s} d\xi\\
&=\int_{\R}f(x) \calF^{-1}\left(\widehat{\dxi }(1+|\xi|^2)^{s}\right),
\end{align*}
with $\calF^{-1}$ the inverse Fourier transform. (This uses the fact that the Fourier transform is a bijection on $H_{-s}$.)

Thus  $\delta_{x}  = \calF^{-1}\left(\widehat{\dxi }(1+|\xi|^2)^{s}\right)\in H_{-s}(\R^n).$  This implies  $\widehat{\dxi } = 
 \calF(\delta_{x}) (1+|\xi|^2)^{-s}$. Using $\calF^{-1}(\calF( f) g) = f* \calF^{-1}(g)$, where $*$ is convolution, we get
\begin{align*}\lefteqn{\dxi (x)}\\
&= \left(\delta_{x} * \calF^{-1}\left((1+|\xi|^2)^{-s}\right)\right)(x) =
\lim_{\sigma\to 0} \int_{\R^n} N(x,\sigma)(y)\calF^{-1}\left((1+|\xi|^2)^{-s}\right)(x-y) dy\\
&= \calF^{-1}\left((1+|\xi|^2)^{-s}\right)(x-x) = 
\frac{1}{(2\pi)^n}\int_{\R^n} e^{i(x-x)\cdot \xi}(1+|\xi|^2)^{-s} d\xi\\
&= \frac{1}{(2\pi)^n}\int_{\R^n} \cos((x-x)\cdot \xi)(1+|\xi|^2)^{-s} d\xi,
%= \frac{1}{\pi}\int_0^\infty \cos(x\xi)(1+|\xi|^2)^{-s/2} d\xi,
\end{align*}
%Set $g(x) = \calF^{-1}\left((1+|\xi|^2)^{-s/2}\right).$
%Then
%$\dxi (x) \approx \int_\R N(x,\sigma)(y)g(x-y)dy \approx g(x-x),$ with
%$$g(x) = \frac{1}{2\pi}\int_\R e^{ix\xi}(1+\xi^2)^{-s/2} d\xi = \frac{1}{2\pi}\int_{\R} \cos(x\xi)(1+\xi^2)^{-s/2} d\xi
%= \frac{1}{\pi}\int_0^\infty \cos(x\xi)(1+\xi^2)^{-s/2} d\xi,
%$$
since $\sin$ is an odd function.  
%Note that while $\delta_{x}$ is considered to have support at $x$, $\dxi$ has support on all of $\R.$

We begin with the case $n=1$, so we are considering $H_s(\R).$ 
%Set $A = s/2.$  We assume
% $s\in 2\mathbb Z$ for now. (Below, in  Case II, we will have $s\in 2\Z + 1$, but then we set
%$A = (s-n+1)/2$, so again $A\in \Z.$  Thus we can apply Lemma \ref{lemma:one}.)  

\begin{lemma} \label{lemma:one}
$$
d_x(y) = \frac{1}{2\pi}
\int_{-\infty}^\infty e^{i|x-y|\xi} (1+\xi^2)^{-s} d\xi = 
%\frac{\pi e^{-|x-y|} |x-y|^{A-1}}{2^{A-1}(A-1)!}
\frac{ e^{-|x-y|} |x-y|^{s-1}}{2^{s}(s-1)!}
\sum_{k=0}^{s-1}  \frac{(s+k-1)!}{k! (s-k-1)!}
 (2|x-y|)^{-k}.$$
 In particular,
 $$ d_x(x) = %2\pi
 \binom{2s-2}{s-1} 2^{-2s+1}.$$
 \end{lemma}
 
We can use $|x-y|$ in place of $x-y$ since $\cos$ is an even function.

\begin{proof}  Fix $a\in \R$.
We compute 
$$\lim_{R\to\infty}\oint_{C_R} e^{iaz} (1+z^2)^{-s},$$
where $C_R$ is the contour given by going along the $x$-axis from $-R$ to $R$ and then the semicircle in the upper half plane centered at the origin and radius $R$, traveled counterclockwise.  By the Jordan Lemma,
the integral over the semicircle goes to zero as $R\to\infty$,
so 
$$\lim_{R\to\infty}\oint_{C_R} e^{iaz} (1+z^2)^{-s} = \int_{-\infty}^\infty e^{iax} (1+x^2)^{-s} dx,$$
the integral we want.

The only pole of the integrand inside the contour is at $z=i$, where $(1+z^2)^{-s} = (i+z)^{-s}(-i +z )^{-s}$ contributes a pole of order $s.$  By the Cauchy residue formula and the standard formula for computing residues, 
\begin{align*} \lefteqn{\oint_{C_R} e^{iaz} (1+z^2)^{-s}} \\
&= 2\pi i \ {\rm Res}_{z=i}  e^{iaz} (1+z^2)^{-s}\\  
%&=  2\pi i \ {\rm Res}_{z=i}  e^{iaz} (z+i)^{-A}\\ 
&= 2\pi i \frac{1}{(s-1)!} \frac {d^{s-1}}{dz^{s-1}}\biggl|_{z=i}(z-i)^s(1+z^2)^{-s}e^{iaz}\\
&= 2\pi i \frac{1}{(s-1)!} \frac {d^{s-1}}{dz^{s-1}}\biggl|_{z=i}(z+i)^{-s}e^{iaz}\\
&=  2\pi i \frac{1}{(s-1)!} \sum_{k=0}^{s-1} \binom{s-1}{k} (-1)^k s(s+1)(s+k-1) (2i)^{-s-k}e^{-a}(ia)^{s-1-k}\\
%\ \ \ \ \ \ \ \ \ \ \ \ \ (*) \\
&=  2\pi ie^{-a} \frac{1}{(s-1)!}\sum_{k=0}^{s-1} (-1)^k \frac{(s-1)!}{k! (s-k-1)!} \frac{(s+k-1)!}{(s-1)!}
i^{-1-2k} 2^{-s-k} a^{s-1-k}\\
%\left(\frac{ia}{2i}\right)^{s-k}\frac{1}{ia}\\
&= \frac{2\pi e^{-a}}{(s-1)!}\sum_{k=0}^{s-1}  \frac{(s+k-1)!}{k! (s-k-1)!}
 2^{-s-k} a^{s-1-k}\\
 &= \frac{\pi e^{-a}a^{s-1}}{2^{s-1}(s-1)!}\sum_{k=0}^{s-1}  \frac{(s+k-1)!}{k! (s-k-1)!}
 2^{-k} a^{-k}.\\
\end{align*}
Letting $R\to\infty$,  replacing $a$ by $|x-y|$ (since $\cos$ is even), and remembering to divide by $2\pi$ %as in (\ref{one}) 
in the definition of the Fourier transform, we get
$$%\int_{-\infty}^\infty e^{i|x-y|\xi} (1+\xi^2)^{-A} d\xi = 
%\frac{\pi e^{-|x-y|}}{2^{A-1}(A-1)!}
\dxix =
\frac{ e^{-|x-y|}|x-y|^{s-1}}{2^{s}(s-1)!}
\sum_{k=0}^{s-1}  \frac{(s+k-1)!}{k! (s-k-1)!}
 (2 |x-y|)^{-k}.$$
 Note that for $x \to y$,  we get a nonzero contribution only for $k = s-1$, so the right hand side equals 
$$%2\pi
\binom{2s-2}{s-1} 2^{-2s+1}.$$
\end{proof}

\begin{remark}  
{\rm The modified Bessel function of the second kind satisfies
$$K_\nu(z) = \frac{\Gamma(\nu + 1/2) (2z)^\nu}{\sqrt{\pi}} \int_0^\infty \cos(at)a^{-2\nu} 
\left(t^2 + \left(\frac{z}{a}\right)^2\right)^{-\nu - 1/2} dt,$$
for $a>0$ and $Re(\nu + 1/2) >0.$ Combining this formula with the half-integer formula
$$K_{(s-1)+ 1/2}(\nu) = \sqrt{\pi/2\nu}\cdot  e^{-\nu} \sum_{k=0}^{s -1}
\frac{(s-1+k)!}{k!(s-1-k)! (2\nu)^{k}}$$
\cite[p.~80, Eq.~(12) and p.~172, Eq.~(1)]{Watson} gives another proof of Lemma 
\ref{lemma:one}. }
\end{remark}

Now we consider the case of $H_s(\R^n,\R)$.
% and assume that $s\in 2\mathbb Z_+.$  
Here
\begin{equation}\label{onea} \dxix = \frac{1}{(2\pi)^n}\int_{\R^n} e^{i  (x-y)\cdot  \xi} 
\left(1+ |\xi|^2\right)^{-s} d\xi.
\end{equation}

Assume $x \neq x.$ Find $B\in SO(n)$ with $Be_n = \frac{x-y}{|x-y|}.$  Then
\begin{align}\label{big}\lefteqn{\int_{\R^n} e^{i  (x-y)\cdot  \xi} 
\left(1+ |\xi|^2\right)^{-s} d\xi} \nonumber\\
&= \int_{\R^n} e^{i  (x-y)\cdot B( \xi)} 
\left(1+ |B(\xi)|^2\right)^{-s} \det (B) d\xi  
=
\int_{\R^n} e^{i  (B^{-1}(x-y))\cdot  \xi} 
\left(1+ |\xi|^2\right)^{-s} d\xi\nonumber\\
& = \int_{\R^n} e^{i  |x-y|  \xi_n} 
\left(1+ |\xi|^2\right)^{-s} d\xi\nonumber\\
&= \int_{\R} e^{i  |x-y|  \xi_n}\left( \int_{\R^{n-1} }  \left(1+\xi_n^2 + \xi_1^2 + \ldots + \xi_{n-1}^2\right)^{-s}
 d\xi_1\ldots d\xi_{n-1}\right) d\xi_n\\
 &= \int_{\R} e^{i  |x-y|  \xi_n}\left( \int_{\R^{n-1} } 
\left(1+\xi_n^2\right)^{-s}  \left(1 +\frac{\xi_1^2}{1+\xi_n^2} + \ldots + \frac{\xi_{n-1}^2}{1+\xi_n^2}\right)^{-s}
 d\xi_1\ldots d\xi_{n-1}\right) d\xi_n\nonumber\\
 &\stackrel{\xi_i\mapsto \xi_i(1+\xi_n^2)^{-1/2}}{=} 
  \int_{\R} e^{i  |x-y|  \xi_n}\left(1 +\xi_n^2\right)^{ -s + (n-1)/2}
  \left( \int_{\R^{n-1} }\left(1+ \xi_1^2 + \ldots + \xi_{n-1}^2 \right)^{-s}
   d\xi_1\ldots d\xi_{n-1}\right) d\xi_n\nonumber\\
   &=\left( \int_{\R} e^{i  |x-y|  \xi_n}\left(1 +\xi_n^2\right)^{ -s + (n-1)/2}d\xi_n \right)
\left(   \int_0^\infty \int_{S^{n-2}} (1+r^2)^{-s} r^{n-2} dr d\theta_1\ldots d\theta_{n-2}\right)\nonumber\\
&=\left( \int_{\R} e^{i  |x-y|  \xi_n}\left(1 +\xi_n^2\right)^{ -s + (n-1)/2}d\xi_n \right)
\left( \int_0^\infty  (1+r^2)^{-s} r^{n-2} dr\right) \vol(S^{n-2})\nonumber\\
&= \left( \int_{\R} e^{i  |x-y|  \xi_n}\left(1 +\xi_n^2\right)^{ -s + (n-1)/2}d\xi_n \right)
\left( \int_0^\infty  (1+r^2)^{-s} r^{n-2} dr\right) 
\frac{2\pi^{\frac{n-1}{2}}}{\Gamma\left(\frac{n-1}{2}\right)}.\nonumber
\end{align}
The first term on the last line is computed in Lemma \ref{lemma:one}, and
is valid for $ -s + (n-1)/2 < -1,$ i.e. for
$$s> (n+1)/2.$$
%\begin{equation}\label{vol} \vol(S^{n-2}) = \frac{2\pi^{\frac{n-1}{2}}}{\Gamma\left(\frac{n-1}{2}\right)},
%\end{equation}
So we have to calculate $\int_0^\infty  (1+r^2)^{-s} r^{n-2} dr$.
\medskip

%\makeatletter
%\newcommand{\pushright}[1]{\ifmeasuring@#1\else\omit\hfill$\displaystyle#1$\fi\ignorespaces}
%\newcommand{\pushleft}[1]{\ifmeasuring@#1\else\omit$\displaystyle#1$\hfill\fi\ignorespaces}
%\makeatother

\noindent {\bf Case I:} $n$ odd %%%%, and we assume $s\in 2\Z.$
\medskip

We use integration by parts repeatedly: % (and assume $s \ge  n+3$):
\begin{align*} \lefteqn{\int_0^\infty  (1+r^2)^{-s} r^{n-2} dr}\\
&= \int_0^\infty \left[ (1+r^2)^{-s}\cdot \frac{2r}{-s + 1}\right] r^{n-3} dr \left(\frac{-s+1}{2}\right)\\
&=- \int_0^\infty (1+r^2)^{-s+1} r^{n-4} dr \cdot \left(\frac{-s+1}{2}\right)(n-3)
\ \ \ \ \ \ \ \ \ \ \ \ \ \ \ \ \ \ \ \ \ \ \ \ \ \ \ \ \ \text{(Step 1)}\\
%&\pushright{\text{(Step 1)}}\\
&=  -\int_0^\infty \left[(1+r^2)^{-s+1}\cdot \frac{2r}{-s+2}\right] r^{n-5} dr \cdot 
\frac{(-s+2)(-s+1)}{2^2}\cdot (n-3)\\
&= \int_0^\infty (1+r^2)^{-s+2} r^{n-6} dr \cdot \frac{(-s+2)(-s+1)}{2^2}\cdot (n-3)(n-5)\ \ \ \ 
\ \ \ \ \  \text{(Step 2)}.
\end{align*}
To get the term $r^1$ in the integrand, we need $\frac{n-3}{2}$ steps. (The exponent drops by $2b+2$ at the $b^{\rm th}$ step, so solve $2b+2 = n-1$ for $b$.)  At Step $\frac{n-3}{2}$, 
the sign in front of the integral is 
$(-1)^{(n-3)/2}$, the exponent of $1+r^2$ is 
$-s + \frac{n-3}{2} = \frac{-2s+n-3}{2}$, and the constant after the integral is
$$\frac{\left(-s + \frac{n-3}{2}\right)\left(-s + \frac{n-5}{2}\right)\cdot\ldots\cdot 
(-s +1)(n-3)(n-5)\cdot\ldots\cdot 2}
{2^{\frac{n-3}{2}}}.$$
(For the final $2$ in the numerator, at the $b^{\rm th}$ step we get $n- (2b+1)= n-(n-3+1) = 2.$)
This equals
%``simplifies" to\footnote{Maybe it's better to recognize the constant as 
%$$(-1)^{(n-3)/2} (s/2-1)!((n-3)/2)!/((s-n+1)/2)! = 
%(-1)^{(n-3)/2} \binom{s/2-1}{(s-n+1)/2}[ ((n-3)/2)!]^2$$. }
\begin{align*}
\lefteqn{ 
\frac{\left(-2s + n-3\right)\left(-2s + n-5\right)\cdot\ldots\cdot 
(-2s+4)(-2s +2)(n-3)!!}{2^{n-3}}  }\\
&= (-1)^{\frac{n-3}{2}} \frac{(2s - n+3)(2s - n+5)\cdot\ldots\cdot 
(2s-4)(2s -2)(n-3)!!}{2^{n-3}}\\
&= (-1)^{\frac{n-3}{2}} \frac{(2s-2)!!(n-3)!!}{2^{n-3} (2s-n+1)!!}.
\end{align*}
Thus
\begin{align*} 
\lefteqn{ \int_0^\infty  (1+r^2)^{-s} r^{n-2} dr}\\
&= (-1)^{\frac{n-3}{2}} \int_0^\infty (1+r^2)^{(-2s+n-3)/2} r dr \cdot 
 (-1)^{\frac{n-3}{2}} \frac{(2s-2)!!(n-3)!!}{2^{n-3} (2s-n+1)!!}\\
 &= \frac{1}{2s-n+1}\cdot \frac{(2s-2)!!(n-3)!!}{2^{n-3} (2s-n+1)!!}.
\end{align*}
%Note that this expression is valid for $s>n.$
Thus for $A = s+ \frac{-n+1}{2}$, with $s \geq (n+1)/2$ an integer,
%\in 2\Z +1, s \geq  n+3$ (e.g., for $s = n+3,$ we have the requirements  $A\in \Z$ and $s>n+1$),
\begin{align}\label{seven}  \lefteqn{\dxix}\\
 &=\frac{1}{(2\pi)^n} 
 \left( \int_{\R} e^{i  |x-y|  \xi}\left(1 +\xi^2\right)^{ -s + (n-1)/2}d\xi_n \right)
 \left(\int_0^\infty  (1+r^2)^{-s} r^{n-2} dr\right)\frac{2\pi^{\frac{n-1}{2}}}{\Gamma\left(\frac{n-1}{2}\right)}  \nonumber\\
 &= \frac{1}{(2\pi)^n} 
 \left( \int_{\R} e^{i  |x-y|  \xi}\left(1 +\xi^2\right)^{ -s + (n-1)/2}d\xi_n \right)
 \cdot \frac{1}{2s-n+1}\cdot \frac{(2s-2)!!(n-3)!!}{2^{n-3} (2s-n+1)!!}
 \cdot  \frac{2\pi^{\frac{n-1}{2}}}{\Gamma\left(\frac{n-1}{2}\right)}\nonumber\\
 &= \frac{1}{\pi^{\frac{n+1}{2}} }\left( \int_{\R} e^{i  |x-y|  \xi}\left(1 +\xi^2\right)^{ -s + (n-1)/2}d\xi \right)
 \cdot \frac{1}{2s-n+1}\cdot \frac{(2s-2)!!(n-3)!!}{2^{2n-4} (2s-n+1)!!}
 \cdot  \frac{1}%2\pi^{\frac{-n-1}{2}}}
 {\left(\frac{n-3}{2}\right)!} \nonumber\\
 &=\frac{1}{\pi^{\frac{n+1}{2}} } \frac{ e^{-|x-y|} |x-y|^{A-1}}{2^{A}(A-1)!}
\sum_{k=0}^{A-1}  \frac{(A+k-1)!}{k! (A-k-1)!}
 (2|x-y|)^{-k}    \nonumber\\
 &\qquad  \cdot \frac{1}{2s-n+1}\cdot \frac{(2s-2)!!(n-3)!!}{2^{2n-4} (2s-n+1)!!}
 \cdot  \frac{1} {\left(\frac{n-3}{2}\right)!}.\nonumber
 \end{align}

\medskip
To treat $x=x$,  we can  use the fact that $\dxix$ is highly differentiable, and let $x\to y$ in the last formula.  As in the one-dimensional case, we get
$$d_x(y) =
\frac{1}{\pi^{\frac{n+1}{2}} } \binom{2s-n-1}{\frac{2s-n-1}{2}}2^{-2s+n}
\frac{1}{2s-n+1}\cdot \frac{(2s-2)!!(n-3)!!}{2^{2n-4} (2s-n+1)!!}
 \cdot  \frac{1} {\left(\frac{n-3}{2}\right)!}.$$

%\begin{align*} d_x(y) &= \frac{1}{(2\pi)^n}\int_{\R^n} \left(1+ |\xi|^2\right)^{-s/2} d\xi = 
%\frac{1}{(2\pi)^n}\int_0^\infty (1+r^2)^{-s/2} r^{n-1} dr d\theta_1\ldots d\theta_{n-1}\\
%&= \frac{1}{(2\pi)^n}\left(\int_0^\infty (1+r^2)^{-s/2} r^{n-1} dr \right) \vol(S^{n-1}).
%\end{align*}

\medskip
\noindent {\bf Case II:} $n$ even %and we assume $s\in 2\Z +1$.
\medskip

Now it takes $\frac{n}{2} - 1$ steps to reduce $\int_0^\infty (1+r^2)^{-s} r^{n-2} dr$ to a constant times\\
$\int_0^\infty (1+r^2)^{-s + n/2 -1 } dr.$ This gives
\begin{align*} 
\lefteqn{ \int_0^\infty  (1+r^2)^{-s} r^{n-2} dr}\\
&= (-1)^{\frac{n}{2}-1} \int_0^\infty (1+r^2)^{-s+n/2-1}  dr\\
&\qquad \cdot 
\frac{\left(-s + \frac{n}{2}-1 \right)\left(-s + \frac{n}{2} - 2\right)\cdot\ldots\cdot 
(-s +1)(n-3)(n-5)\cdot\ldots\cdot 1}
{2^{\frac{n}{2}-1}}\\
&= \int_0^\infty (1+r^2)^{-s+n/2-1}  dr \cdot \frac{(s-1)(s-3)\cdot\ldots\cdot(s-\frac{n}{2}+1)
\left(n-3\right)!!}{2^{\frac{n}{2}-1}}\\
%\left(\frac{n-3}{2}\right)!!}{2^{\frac{n}{2}-1}}\\
%! \left(\frac{n-3}{2}\right)!}{\left(\frac{s-n}{2}\right)!}\\
%&= \pi \binom{s-n}{\frac{s-n}{2}} 2^{n-s} \frac{(s/2-1)! \left(\frac{n-3}{2}\right)!}{\left(\frac{s-n}{2}\right)!},
\end{align*}
%by Lemma \ref{lemma:one} for the case $x = x$ with $A = (s-n+2)/2.$\footnote{Here  $((n-3)/2)!$ means 
%$\Gamma((n-1)/2).$  Maybe it's best to use $(n-3)!!$ in the even case.}
By %the substitution $r = \tan(u)$, 
Wolfram Alpha, we get
%\footnote{Use Wolfram Alpha for  ``int cos\^\ \{2k+1\}(x) dx" and 
%``Gauss hypergeometric functions at (1/2, k+1, k+2, 1).''}
$$\int_0^\infty (1+r^2)^{-k} = %\int_0^{\pi/2} \cos^{2k-1}(u) du
%= {}_2F_1(1/2,k+1; k+2, 1) 
 \frac{\sqrt{\pi}\ \Gamma(k- \frac{1}{2})}%{\Gamma\left(k+\frac{3}{2}\right)}
{2\Gamma(k)},$$
so
$$ \int_0^\infty  (1+r^2)^{-s/2} r^{n-2} dr = \frac{\sqrt{\pi}\Gamma(s-\frac{n+1}{2})}{\Gamma(s-\frac{n}{2}+1)}
\cdot\frac{(s-1)(s-3)\cdot\ldots\cdot(s-\frac{n}{2}+1)
\left(n-3\right)!!}{2^{\frac{n}{2}-1}}.
$$
%%%%%%%%%%%%%%%%%%%%%%%%%%%%%%%%%%%%%%
Using $\Gamma(n+ (1/2)) = \frac{(2n)!}{4^n n!} \sqrt{\pi}$, we obtain
%Therefore,\footnote{Using $\Gamma(n+ (1/2)) = \frac{(2n)!}{4^n n!} \sqrt{\pi}.$}
 for $k= (2s-n+3)/2$,\begin{align*} 
\lefteqn{ \int_0^\infty  (1+r^2)^{-s/2} r^{n-2} dr}\\
&=   \int_0^\infty (1+r^2)^{(-s+n-2)/2}  dr \cdot 
\frac{(s/2-1)(s/2-1)\cdot\ldots\cdot(s/2-\frac{n}{2}+1)
\left(\frac{n-3}{2}\right)!!}{2^{\frac{n}{2}-1}}\\
%(s/2-1)(s/2-1)\cdot\ldots\cdot(s/2-\frac{n}{2}+1)
%\cdot\left(\frac{n-3}{2}\right)!\\
&=  \frac{\sqrt{\pi}\ \Gamma\left(\frac{s-n+5}{2}\right)}{\Gamma\left(\frac{s-n+4}{2}\right)}
\frac{(s/2-1)(s/2-1)\cdot\ldots\cdot(s/2-\frac{n}{2}+1)
\left(\frac{n-3}{2}\right)!!}{2^{\frac{n}{2}-1}}\\
% (s/2-1)(s/2-1)\cdot\ldots\cdot(s/2-\frac{n}{2}+1)\cdot\left(\frac{n-3}{2}\right)!\\
&=  2^{s-n+3}\left[\binom{s-n+3}{\frac{s-n+3}{2}} \right]^{-1}\cdot
\frac{(s/2-1)(s/2-1)\cdot\ldots\cdot(s/2-\frac{n}{2}+1)
\left(\frac{n-3}{2}\right)!!}{2^{\frac{n}{2}-1}}.
% (s/2-1)(s/2-1)\cdot\ldots\cdot(s/2-\frac{n}{2}+1)\cdot\left(\frac{n-3}{2}\right)!
\end{align*}

By (\ref{big}), we have for $A = (s-n+1)/2$ and $x\neq x$, 
\begin{align}\label{eight}
\lefteqn{\dxix =\frac{1}{(2\pi)^n}\int_{\R^n} e^{i  (x-x)\cdot  \xi} 
\left(1+ |\xi|^2\right)^{-s/2} d\xi} \nonumber\\
&= \frac{1}{(2\pi)^n}\left( \int_{\R} e^{i  |x-y|  \xi_n}\left(1 +\xi_n^2\right)^{ -s/2 + (n-1)/2}d\xi_n \right)
\left( \int_0^\infty  (1+r^2)^{-s/2} r^{n-2} dr\right) \vol(S^{n-2})\nonumber \\
&= \frac{1}{(2\pi)^n}\left(\frac{ e^{-|x-y|} |x-y|^{A-1}}{2^{A}(A-1)!}
\sum_{k=0}^{A-1}  \frac{(A+k-1)!}{k! (A-k-1)!}
 (2|x-y|)^{-k}\right)\\
 &\qquad \cdot\left(  2^{s-n+3}\left[\binom{s-n+3}{\frac{s-n+3}{2}} \right]^{-1}\cdot
\frac{(s/2-1)(s/2-1)\cdot\ldots\cdot(s/2-\frac{n}{2}+1)
\left(\frac{n-3}{2}\right)!!}{2^{\frac{n}{2}-1}} \right)\nonumber \\
 &\qquad 
 %\cdot\left( 2^{s-n+3}\left[\binom{s-n+3}{\frac{s-n+3}{2}} \right]^{-1}\cdot
% (s/2-1)(s/2-1)\cdot\ldots\cdot(s/2-\frac{n}{2}+1)
%\cdot\left(\frac{n-3}{2}\right)!\right)
% \pi \binom{s-n}{\frac{s-n}{2}} 2^{n-s} \frac{(s/2-1)! \left(\frac{n-3}{2}\right)!}{\left(\frac{s-n}{2}\right)!} \right) 
\cdot \frac{2\pi^{\frac{n-1}{2}}}{\Gamma\left(\frac{n-1}{2}\right)}\nonumber\\
&= \frac{2^{(s+5)/2} \left[\left(\frac{s-n+3}{2}\right)!\right]^2 (s/2-1)(s/2-3)\cdot...\cdot(s/2 - n/2 +1)(n-1)!}{\pi^{(n/2)+1} (s-n+3)! \left( \frac{s-n-1}{2}\right)! (2n-2)!}\nonumber\\
&\qquad  \cdot
 e^{-|x-y|} |x-y|^{A-1}
\sum_{k=0}^{A-1}  \frac{(A+k-1)!}{k! (A-k-1)!}
 (2|x-y|)^{-k} \nonumber\\
 & := C_{s,n}e^{-|x-y|} |x-y|^{A-1}
\sum_{k=0}^{A-1}  \frac{(A+k-1)!}{k! (A-k-1)!}
 (2|x-y|)^{-k}. \nonumber
\end{align}
 Again, to have $n$ even, $s\in 2\Z +1$, and $s > n+1$, we need $s \geq n+3$.

For $x=x$, we get for the only nonzero term $k = A-1 = (s-n-1)/2$,
\begin{align*} d_x(y) &= C_{s,n} \frac{(s-n-1)!}{\left(\frac{s-n+1}{2}\right)!} 2^{(-s+n+1)/2},
\end{align*}
which can be somewhat simplified.

\vskip 0.1in
\noindent {\bf Simplifying the results} 
\vskip 0.1in

The basic fact is that $H_{n+k +\epsilon} \subset C^{ n/2 + k}(\R^n)$  for any $\epsilon >0.$ 
Recall that we can choose $s = n+3$, in which case $\dxix\in 
H_{n+3} \subset C^{[n/2] +1}(\R^n).$\\

\noindent {\bf Case I:} $n$ odd, $s = n+3$, $A = (s-n+1)/2 = 2$\\

By (\ref{seven}), we get
\begin{align}\label{nine}
\dxix &=  C_n e^{-|x-y|} |x-y| \left(1 + |x-y|^{-1}\right) = C_n e^{-|x-y|}\left(1+  |x-y|\right);\\
C_n &=  
\frac{1}{\pi^{\frac{n+1}{2}} } \frac{ 1}{2^{2}}
\frac{1}{4}\cdot \frac{(n+1)!!(n-3)!!}{2^{2n-4} (4)!!}
 \cdot  \frac{1} {\left(\frac{n-3}{2}\right)!}.
\nonumber
\end{align}
The dimension constant $C_n$ can be simplified, since
\begin{align*}(n-3)!! &= 2^{(n-3)/2}\frac{n-3}{2}\cdot \frac{n-3}{5}\cdot...\cdot 1\\
 \Rightarrow
C_n &= \frac{1}{\pi^{\frac{n+1}{2}} }\frac{1}{2^7} \frac{2^{(n+1)/2}2^{(n-3)/2} (n-1)!}{2^{2n-4}}
=  \frac{1}{\pi^{\frac{n+1}{2}} }\frac{(n-1)!}{2^{n+4}}.  
%\ \ \ \textcolor{red}{\text{ [check\ this!]} }
\end{align*}
In any case, $\dxix$ for $H_{n+3}(\R^n)$ is a simple expression times a dimension constant. 
\medskip

\noindent {\bf Case II:} $n$ even, $s = n+3$, $A = (s-n+1)/2 = 2$\\

By (\ref{eight}), we get
\begin{align}\label{ten}
\dxix &= \frac{2^{(n/2)+5}(3!)^2 \left(\frac{n+1}{2}\right) \left( \frac{n-3}{2}\right)\left( \frac{n-7}{2}\right)
\cdot ...\cdot \left(\frac{5}{2}\right) (n-1)!}{\pi^{(n/2)+1} 6! (2n-2)!}
\left(e^{-|x-y|}\left(1+  |x-y|\right) \right).
\end{align}

\bibliographystyle{plain}
\bibliography{iclr2023_conference}

\end{document}

%% file: math_commands.tex
\def\eqref#1{equation~\ref{#1}}
\def\1{\bm{1}}

\DeclareMathAlphabet{\mathsfit}{\encodingdefault}{\sfdefault}{m}{sl}
\SetMathAlphabet{\mathsfit}{bold}{\encodingdefault}{\sfdefault}{bx}{n}

\newcommand{\R}{\mathbb{R}}